\documentclass{amsart}
\usepackage{amssymb,amsmath,latexsym,epsfig}

\newtheorem{thm}{Theorem}[section]

\newtheorem{cor}[thm]{Corollary}

\newtheorem{lemma}[thm]{Lemma}
\newtheorem{prob}[thm]{Problem}

\newcommand{\del}{\backslash}
\newcommand{\cl}{\hbox{\rm cl}}
\newcommand{\GF}{\hbox{\rm GF}}
\newcommand{\PG}{\hbox{\rm PG}}

\newcommand{\real}{\mathbb{R}}
\newcommand{\cw}[1]{\mathrm{CW}(#1)}
\newcommand{\mc}[1]{\mathcal{#1}}
\newcommand{\join}{\lor}
\newcommand{\meet}{\land}
\newcommand{\frp}{\mathbin{\Box}}

\title[The Lattice of Cyclic Flats]{The Lattice of Cyclic Flats of a Matroid}
\date{\today}
\author[Joseph E.~Bonin]
       {Joseph E.~Bonin}
\address[Joseph E.~Bonin]
{Department of Mathematics\\ The George Washington University\\
Washington, D.C. 20052, USA} \email [Joseph E.~Bonin] {jbonin@gwu.edu}
\author[Anna de Mier]
{Anna de Mier}
\address[Anna de Mier]
{Mathematical Institute\\  24--29 St~Giles' \\ 
Oxford OX1 3LB, United Kingdom} \email [Anna de Mier] {ademier@gmail.com}

\begin{document}

\begin{abstract}
  A flat of a matroid is cyclic if it is a union of circuits.  The
  cyclic flats of a matroid form a lattice under inclusion.  We study
  these lattices and explore matroids from the perspective of cyclic
  flats.  In particular, we show that every lattice is isomorphic to
  the lattice of cyclic flats of a matroid.  We give a necessary and
  sufficient condition for a lattice $\mc{Z}$ of sets and a function
  $r:\mc{Z}\rightarrow \mathbb{Z}$ to be the lattice of cyclic flats
  of a matroid and the restriction of the corresponding rank function
  to $\mc{Z}$.  We apply this perspective to give an alternative view
  of the free product of matroids and we show how to compute the Tutte
  polynomial of the free product in terms of the Tutte polynomials of
  the constituent matroids.  We define cyclic width and show that this
  concept gives rise to minor-closed, dual-closed classes of matroids,
  two of which contain only transversal matroids.
\end{abstract}

\maketitle


\section{Introduction}\label{sec:intro}

\begin{figure}[t]
\begin{center}
\epsfxsize 10cm \epsffile{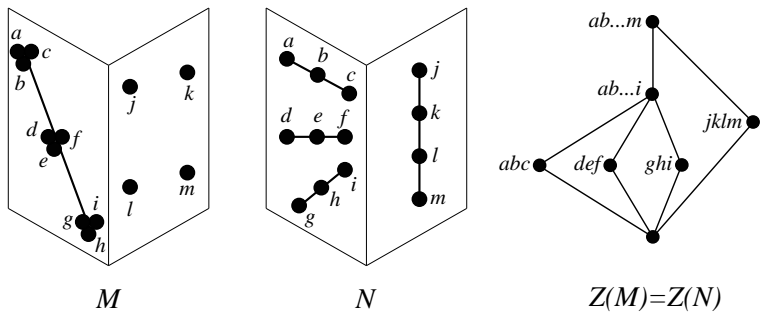} 
\caption{Two matroids that have the same lattice 
  of cyclic flats.}\label{fig:lattice}
\end{center}
\end{figure}

A flat of a matroid is \emph{cyclic} if it is a (possibly empty) union
of circuits.  Cyclic flats have played several important roles in
matroid theory, starting with the theory of transversal matroids (see,
for example, \cite{affine,ingleton}).  The cyclic flats of a matroid
$M$, ordered by inclusion, form a lattice, $\mc{Z}(M)$: the join of
two cyclic flats $A$ and $B$ in $\mc{Z}(M)$ is $\cl(A\cup B)$ and
their meet is the union of all circuits contained in $A\cap B$.  (See
Figure~\ref{fig:lattice}.)  This paper studies the lattice of cyclic
flats and offers insights into some topics in matroid theory from this
perspective.

Figure~\ref{fig:lattice} shows that properties of the lattice of
flats, such as being graded, may fail in the lattice of cyclic flats.
Indeed, in~Section~\ref{sec:basic} we prove that every lattice is
isomorphic to the lattice of cyclic flats of some matroid.
Figure~\ref{fig:lattice} also shows that the collection of cyclic
flats alone does not determine the matroid.  However, the cyclic flats
together with their ranks determine the
matroid~\cite[Proposition~2.1]{affine}.  Section~\ref{sec:axioms}
gives a necessary and sufficient condition for a collection $\mc{Z}$
of sets and a function $r:\mc{Z}\rightarrow \mathbb{Z}$ to be the
lattice of cyclic flats of a matroid and the restriction of the
corresponding rank function to $\mc{Z}$.  This axiom scheme for
matroids is used in Section~\ref{sec:free} to give a simple,
alternative perspective on the free product~\cite{free,free2}.
Section~\ref{sec:cw} introduces the concept of cyclic width and shows
that the class of matroids of cyclic width $k$ or less is closed under
minors and duals. We explore these classes for small values of $k$.
The final section contains open problems.

We assume familiarity with matroid theory.  We follow the notation
of~\cite{oxley}.  We call a matroid \emph{bitransversal} if it is both
transversal and cotransversal (or a strict gammoid).  Recall that the
free extension $M+e$ of $M$ by $e$, for $e\not\in E(M)$, is the
matroid on $E(M)\cup e$ whose circuits are the circuits of $M$
together with the sets $B\cup e$ for bases $B$ of $M$.  Free
coextension is the dual operation, $(M^*+e)^*$.

Our lattice theory notation follows~\cite{hilary} except that we use
$0_{\mc{Z}}$ and $1_{\mc{Z}}$ for the least and greatest elements of a
lattice $\mc{Z}$.  As is true of a lattice of flats and a lattice of
cyclic flats, the elements of most lattices in this paper are sets and
the order relation is containment; we use the term \emph{lattice of
  sets} for such a lattice.  The meet and join operations of a lattice
of sets need not be intersection and union.  All lattices considered
in this paper are finite.


We close this introduction by recalling that cyclic flats link two
notions of duality.  Since cyclic flats are unions of circuits as well
as intersections of hyperplanes, $X$ is a cyclic flat of $M$ if and
only if $E(M)-X$ is a cyclic flat of $M^*$.  Thus, as noted
in~\cite{affine,ingleton}, $\mc{Z}(M^*)$ is isomorphic to the dual (as
a lattice) of $\mc{Z}(M)$.


\section{The Lattice of Cyclic Flats}\label{sec:basic}


By considering matroids in which every element is in a nontrivial
parallel class, it follows that every geometric lattice is isomorphic
to the lattice of cyclic flats of some matroid.  The main result of
this section is that for any lattice $\mc{Z}$, there is a matroid $M$
(in fact, a bitransversal matroid) for which the lattice $\mc{Z}(M)$ of
cyclic flats is isomorphic to $\mc{Z}$.  We also show that the counterpart of
this result does not hold for matroids that are representable over a
fixed finite field.

We will use fundamental transversal matroids and their simplex
representation (\cite{affine,ingleton} or~\cite[Section~12.2]{oxley}).
A \emph{fundamental} (or \emph{principal}) \emph{transversal matroid}
is a matroid $M$ that has a basis $B$ such that each cyclic flat of
$M$ is spanned by a subset of $B$.  To construct a geometric
representation of $M$, let $B$ be $\{v_1,v_2,\ldots,v_n\}$ and put
$v_1,v_2,\ldots,v_n$ at the $n$ vertices of an $n$-simplex in
$\real^n$.  For each $x$ in $E(M)-B$, place $x$ freely in the
(possibly empty) face of the simplex spanned by $C(x,B)-x$, where
$C(x,B)$ is the fundamental circuit of $x$ with respect to $B$.  Note
that if $X$ is a rank-$k$ cyclic flat of $M$, then $X$ consists of the
points in a $(k-1)$-dimensional face of the simplex.  Conversely, any
matroid with a geometric representation on a simplex that has elements
at each vertex and for which each rank-$k$ cyclic flat consists of the
points in a $(k-1)$-dimensional face of the simplex is a fundamental
transversal matroid.  Fundamental transversal matroids are indeed
transversal matroids; in fact, every transversal matroid is a deletion
of a fundamental transversal matroid. It is well known that these
matroids are bitransversal.

We are ready to treat the main result of this section.

\begin{thm}\label{thm:alllattices}
  Every lattice is isomorphic to the lattice of cyclic flats of a
  bitransversal matroid.
\end{thm}

\begin{proof}
  Let $\mc{Z}$ be a lattice.  We construct a fundamental transversal
  matroid $M$ for which $\mc{Z}(M)$ is isomorphic to $\mc{Z}$.  We
  first realize $\mc{Z}$ as a lattice of sets. Let $B$ be the set of
  elements in $\mc{Z}$ other than $1_{\mc{Z}}$.  For $z$ in $\mc{Z}$,
  let $V_{z}$ be $\{y: y\not\geq z\}$.  Thus, $V_z\subseteq B$.  The
  map $z\mapsto V_z$ is an isomorphism between $\mc{Z}$ and the
  lattice of sets $\{V_z\,:\, z\in\mc{Z}\}$ since it is a bijection
  and $x\leq z$ in $\mc{Z}$ if and only if $V_x\subseteq V_z$.
  
  Let $M$ be the following fundamental transversal matroid on $B\cup
  \{s_z: z\in \mc{Z}\}$.  Place the elements of $B$ at the vertices of
  a $|B|$-simplex in $\mathbb{R}^{|B|}$ and, for each $z$ in $\mc{Z}$,
  put $s_z$ freely in the face spanned by $V_z$.  By construction,
  $V_z$ spans a cyclic flat of $M$, namely, $V_z\cup \{s_x\,:\, x\leq
  z\}$.  With the isomorphism above, the theorem follows upon showing
  that these are the only cyclic flats of $M$.
  
  Assume the subset $V$ of $B$ spans a cyclic flat $F$ of $M$.  For
  any $s_z$ in $F-V$, we have $C(B,s_z)=V_z\cup s_z\subseteq F$.
  Also, each element $v$ of $V$ is in such a fundamental circuit, for
  if $v$ were not, then all elements of $F-V$ would be in $\cl(V-v)$,
  so $F-v$ would be a flat, contrary to $F$ being cyclic.  Thus, $F$
  is the closure of $\bigcup\,\{V_z\,:\, s_z\in F-V\}$.  Note that
  $V_{x\join z}=V_x\cup V_z$, so $F$ is $\cl(V_{z'})$, where $z'$ is
  $\bigvee\,\{z\,:\,s_z\in F-V\}$.  Thus, as needed, the cyclic flats
  of $M$ are the sets $V_z\cup \{s_x: x\leq z\}$.
\end{proof}

After (re)discovering Theorem~\ref{thm:alllattices}, we learned of the
work of Sims~\cite{julie} extending Dilworth's embedding
theorem~\cite[Theorem 14.1]{cd} to rank-finite independence spaces and
lattices of finite length.  Her approach has much in common with
Dilworth's and gives the lattice of cyclic flats as a sublattice of
the lattice of flats.  Dilworth's proof was very influential in the
development of the theory of submodular functions.  We note that a
minor modification of the proof above gives a short, geometric proof
of a strengthening of Dilworth's theorem: for each finite lattice
$\mc{Z}$, there is a transversal matroid $M$ for which (i) $\mc{Z}$ is
isomorphic to $\mc{Z}(M)$ and (ii) $\mc{Z}(M)$ is a sublattice of the
lattice of flats of $M$.  Indeed, adapt the construction as follows:
for each $z\in \mc{Z}$, instead of putting one point freely in the
flat spanned by $V_z$, put a set $S_z$ of $|V_z|+1$ points freely in
this flat, where $S_z\cap S_x=\emptyset$ for $z\ne x$, and then delete
$B$ to get a transversal matroid $M$.  The cyclic flats of $M$ are the
sets $F_z=\bigcup_{y\leq z}S_y$; also, $F_x\cap F_y=F_{x\meet y}$, so
the meet in $\mc{Z}(M)$ is that of the lattice of flats; the joins in
these lattices always agree, thus giving the sublattice assertion.

Since transversal matroids are representable over the reals, the
counterpart of Theorem~\ref{thm:alllattices} holds for real matroids.
It is natural to ask what other classes of matroids have a counterpart
of this result.  Our next goal (Theorem~\ref{thm:gfqdonotcover}) is to
show that matroids representable over $\GF(q)$ do not have this
property.

The proof of Theorem~\ref{thm:gfqdonotcover} uses the class of nested
matroids, which also appears in Section~\ref{sec:cw}.  A \emph{nested
  matroid} is a matroid that can be obtained from the empty matroid by
iterating the operations of adding isthmuses and taking free
extensions. Nested matroids have appeared several times in the
literature (see~\cite[Section 4]{lpm2} for references). They are a
subclass of bitransversal matroids.

The following observations about the operations used to construct
nested matroids will be useful.  For a single-element extension $M'$
of $M$, the lattices $\mc{Z}(M)$ and $\mc{Z}(M')$ are equal if and
only if the element of $E(M')-E(M)$ is an isthmus of $M'$.  The
lattice $\mc{Z}(M+e)$ of the free extension of $M$ is formed from
$\mc{Z}(M)$ by removing $E(M)$, if it is in $\mc{Z}(M)$, and adjoining
$E(M)\cup e$ to $\mc{Z}(M)$.  If $\mc{Z}$ is the lattice of cyclic
flats of some matroid and $T$ is disjoint from $1_{\mc{Z}}$, then the
following construction gives all matroids $M$ for which $\mc{Z}(M)$ is
$\mc{Z}\cup\{1_{\mc{Z}}\cup T\}$:
\begin{itemize}
\item[(1)] start with any matroid $M_0$ for which $E(M_0)$ is
  $1_{\mc{Z}}$ and $\mc{Z}(M_0)$ is $\mc{Z}$,
\item[(2)] partition $T$ into two nonempty subsets $T_i$ and $T_f$,
\item[(3)] take the direct sum of $M_0$ and the free matroid on $T_i$
  to get $M_1$,
\item[(4)] take free extensions of $M_1$ by the elements in $T_f$ to
  get $M_2$, and finally
\item[(5)] let $M$ be the direct sum of $M_2$ and any free matroid.
\end{itemize}
Note that in steps (2)--(4), only the cardinalities $|T_i|$ and
$|T_f|$, not the elements in these sets, matter.  These observations
give us the following result, which is essentially Lemma~2
of~\cite{opr}. Recall that a chain is a linearly ordered set.

\begin{lemma}\label{lem:nested}
  A matroid $M$ is nested if and only if $\mc{Z}(M)$ is a chain.
\end{lemma}

We now show that some lattices do not arise as the lattice of cyclic
flats of any matroid that is representable over $\GF(q)$.

\begin{thm}\label{thm:gfqdonotcover}
  For a prime power $q$, if $\mc{Z}(M)$ is a chain with $q+2$ or more
  elements, then $M$ is not representable over $\GF(q)$.
\end{thm}

\begin{proof}
  We show that if $\mc{Z}(M)$ is a chain with $k+2$ elements, then $M$
  has $U_{k,k+2}$ as a minor; thus if $k\geq q$, then $M$ is not
  representable over $\GF(q)$.  Let the cyclic flats of $M$ be
  $X_0\subset X_1\subset \cdots \subset X_{k+1}$.  By
  Lemma~\ref{lem:nested}, $M$ is a nested matroid, so $M|X_j$, for
  $1\leq j\leq k+1$, is obtained by partitioning $X_j-X_{j-1}$ into
  nonempty sets $I_j$ and $F_j$, taking the direct sum of $M|X_{j-1}$
  and the free matroid on $I_j$, and then taking free extensions by
  the elements in $F_j$.  Therefore the minor
  $$M|X_{k+1}/I_{k+1}\del (X_0\cup F_1\cup F_2\cup \cdots \cup
  F_{k-1})$$
  has just two cyclic flats, namely, $\emptyset$ and
  $I_1\cup I_2\cup\ldots \cup I_k\cup F_k\cup F_{k+1}$, and so is a
  uniform matroid; that the rank and nullity of this uniform minor are
  at least $k$ and $2$, respectively, completes the proof.
\end{proof}


\section{An Axiom Scheme for Cyclic Flats}\label{sec:axioms}


In this section we formulate the definition of a matroid using cyclic
flats and their ranks.  The mixed character of this axiom scheme sets
it apart from most others: it is built jointly on lattice properties
(but using only some of the flats) and partial information about the
rank function.

We will use the following characterization of independent sets and
circuits in terms of cyclic flats and their ranks.

\begin{lemma}\label{lem:cyc+rank}
  \emph{(i)} \ A set $I\subseteq E(M)$ is independent in $M$ if and
  only if
  $|I\cap X| \leq r(X)$ for every cyclic flat $X$ of $M$.\\
  \emph{(ii)} \ A set $C\subseteq E(M)$ is a circuit of $M$ if and
  only if $C$ is minimal with the property that there is a cyclic flat
  $X$ with $C\subseteq X$ and $|C|=r(X)+1$.
\end{lemma}

\begin{proof}
  Statement (i) is immediate.  As the minimal dependent sets, the
  circuits are the minimal sets $C$ with $|C\cap X|>r(X)$ for some
  cyclic flat $X$; since $X$ can be $\cl(C)$, this property can be
  replaced by the simpler one in statement (ii).
\end{proof}

To motivate our main result, we note three properties of the lattice
$\mc{Z}(M)$ and the rank function $r$ of $M$.  First, $0_{\mc{Z}(M)}$,
which is $\cl(\emptyset)$, has rank $0$.  Next, if $X$ and $Y$ are in
$\mc{Z}(M)$ and $X\subsetneq Y$, then $0<r(Y)-r(X)<|Y-X|$ since $M|Y$
has no isthmuses.  Finally, since for cyclic flats $X$ and $Y$ we have
$r(X\join Y)=r(X\cup Y)$ and $r(X\cap Y)=r(X\meet Y) +|(X\cap Y) -
(X\meet Y)|$, we get the following specialization of semimodularity:
$$r(X)+r(Y)\geq r(X\join Y) + r(X\meet Y) + |(X\cap Y) - (X\meet
Y)|.$$
We next show that these properties give an axiom scheme for
matroids.

\begin{thm}\label{thm:axioms}
  Let $\mc{Z}$ be a collection of subsets of a set $S$ and let $r$ be
  an integer-valued function on $\mc{Z}$.  There is a matroid for
  which $\mc{Z}$ is the collection of cyclic flats and $r$ is the rank
  function restricted to the sets in $\mc{Z}$ if and only if
\begin{itemize}
\item[(Z0)] $\mc{Z}$ is a lattice under inclusion,
\item[(Z1)] $r(0_{\mc{Z}})=0$,
\item[(Z2)] $0<r(Y)-r(X)<|Y-X|$ for all sets $X,Y$ in $\mc{Z}$ with
  $X\subsetneq Y$, and  
\item[(Z3)] for all sets $X,Y$ in $\mc{Z}$,
  $$r(X)+r(Y)\geq r(X\join Y) + r(X\meet Y) + |(X\cap Y) - (X\meet
  Y)|.$$
\end{itemize}
\end{thm}

\begin{proof}
  We have seen that these properties are necessary; we focus on
  sufficiency.  Assume $\mc{Z}$ and $r$ satisfy properties (Z0)--(Z3).
  We define a collection $\mc{C}$ (following Lemma~\ref{lem:cyc+rank})
  that we show is the set of circuits of a matroid on $S$.  We then
  show that the cyclic flats of this matroid are the sets in $\mc{Z}$
  and that $r$ gives their ranks.


\begin{quote}
  {\bf (\ref{thm:axioms}.1)} The collection $\mc{C}$ of all minimal
  subsets $C$ of $S$ for which there is a set $X\in \mc{Z}$ with
  $C\subseteq X$ and $|C|=r(X)+1$ is the collection of circuits of a
  matroid.
\end{quote}

\begin{proof}[Proof of (\ref{thm:axioms}.1).]
  We focus on circuit elimination since the other two circuit axioms
  clearly hold.  Thus, let $C$ and $C'$ be distinct sets in $\mc{C}$
  with $a$ in $C\cap C'$; let $X$ and $X'$ in $\mc{Z}$ contain $C$ and
  $C'$, respectively, with $|C|=r(X)+1$ and $|C'|=r(X')+1$.  To show
  that $(C\cup C')-a$ contains a set in $\mc{C}$, it suffices to prove
  the inequality
  $$
  |(C\cup C')-a|>r(X\join X').
  $$
  Property (Z3) and the equalities $|C|=r(X)+1$ and $|C'|=r(X')+1$
  give
  \begin{align}
    |(C\cup C')-a| = &\,\, |C| + |C'| - |C\cap C'| - 1 \notag \\
    > &\,\, r(X) + r(X') - |C\cap C'|\notag \\
    \geq &\,\, r(X\join X') + r(X\meet X') + |(X\cap X')-(X\meet X')|
    - |C\cap C'|.\notag
  \end{align}
  Thus, the desired inequality follows upon showing that the sum of
  the last three terms is nonnegative, that is, $|C\cap C'|\leq
  r(X\meet X') + |(X\cap X')-(X\meet X')|$.  Now $|(C\cap C')\cap
  (X\meet X')|\leq r(X\meet X')$ since $X\meet X'$ is in $\mc{Z}$ and
  no subset of $C\cap C'$ is in $\mc{C}$, so 
  \begin{align}
    |C\cap C'| = &\,\, |(C\cap C')\cap (X\meet X')| + |(C\cap C') -
    (X\meet X')|
    \notag\\
    \leq &\,\, r(X\meet X') + |(X\cap X') - (X\meet X')|,\notag
  \end{align}
  and therefore circuit elimination holds.
\end{proof}

Let $M_{\mc{Z}}$ be the matroid on $S$ that has $\mc{C}$ as its
collection of circuits and let $r_{\mc{Z}}$ be its rank function.  The
next step in essence identifies some elements of $\mc{Z}$ as the
closures of circuits.  


\begin{quote}
  {\bf (\ref{thm:axioms}.2)} For any circuit $C$ of $M_{\mc{Z}}$,
  there is a unique set $X\in \mc{Z}$ with $C\subseteq X$ and
  $r(X)=|C|-1$.  Furthermore, if $Y\in \mc{Z}$ and $C\subseteq Y$,
  then $X\subseteq Y$.
\end{quote}

\begin{proof}[Proof of (\ref{thm:axioms}.2).]
  By the definition of $\mc{C}$, there is a set $X$ in $\mc{Z}$ with
  $C\subseteq X$ and $r(X)=|C|-1$. Uniqueness follows by showing that
  if $Y\in \mc{Z}$ and $C\subseteq Y$, then $X\subseteq Y$.  Assume
  instead $X\not\subseteq Y$.  Thus, $X\meet Y\subsetneq X$.  We have
  $C\not\subseteq X\meet Y$, for otherwise, since $|C|>r(X\meet Y)+1$,
  some proper subset of $C$ would be in $\mc{C}$, which is impossible.
  Property~(Z3) and the equality $r(X)=|C|-1$ give
  \begin{equation}\label{ineq}
  r(Y)+|C|-1 \geq r(X\join Y) +r(X\meet Y)+|(X\cap Y)- (X\meet Y)|.
  \end{equation}
  No proper subset of $C$ is in $\mc{C}$, so $r(X\meet Y)\geq |C\cap
  (X\meet Y)|$.  Since $C \subseteq X\cap Y$, we have $|(X\cap Y)-
  (X\meet Y)| \geq |C-(X\meet Y)|$.  Thus, the last two terms in
  inequality~(\ref{ineq}) contribute at least $|C|$, so
  inequality~(\ref{ineq}) implies $r(Y)-1\geq r(X\join Y)$, contrary
  to property~(Z2).  Thus, $X\subseteq Y$.
\end{proof}

We use $\bar{C}$ to denote the unique set in~(\ref{thm:axioms}.2).

The next step is used to prove (\ref{thm:axioms}.4), that $r$ and
$r_{\mc{Z}}$ agree on the sets in $\mc{Z}$.


\begin{quote}
  {\bf (\ref{thm:axioms}.3)} For all $X,Y$ in $\mc{Z}$, the union of
  $(X\cap Y)- (X\meet Y)$ and any basis of $X\meet Y$ is a basis for
  $X\cap Y$, so
  $$r_{\mc{Z}}(X\cap Y)=r_{\mc{Z}}(X\meet Y)+ |(X\cap Y)- (X\meet
  Y)|.$$
\end{quote}

\begin{proof}[Proof of (\ref{thm:axioms}.3).]
  Assume $X$ and $Y$ are incomparable, otherwise the result is
  trivial.  Let $B$ be a basis of $X\meet Y$ and let $B'$ be $B\cup
  \bigl((X\cap Y)- (X\meet Y)\bigr)$. Clearly $B'$ spans $X\cap Y$, so
  it suffices to show that $B'$ is independent. Assume, to the
  contrary, that $B'$ contains a circuit $C$. Thus, $C$ is contained
  in both $X$ and $Y$, so $\bar{C}$ and, consequently, $(X\meet Y)
  \join \bar{C}$ are contained in both $X$ and $Y$.  However, $(X\meet
  Y) \join \bar{C}$ contains elements in $(X\cap Y)-X\meet Y$ since
  $C$ must, so $X\meet Y\subsetneq (X\meet Y) \join \bar{C}$.  This
  contradicts the definition of meet and so proves the claim.
\end{proof}


\begin{quote}
  {\bf (\ref{thm:axioms}.4)} If $X$ is in $\mc{Z}$, then
  $r(X)=r_{\mc{Z}}(X)$.
\end{quote}

\begin{proof}[Proof of (\ref{thm:axioms}.4).]
  By the definition of $\mc{C}$, any subset of $X$ with more than
  $r(X)$ elements contains a circuit, so $r_\mc{Z}(X)\leq r(X)$.
  Thus, it suffices to show that $X$ contains an independent set of
  size $r(X)$.  This statement clearly holds for $0_{\mc{Z}}$.  Let
  $X$ be a minimal set in $\mc{Z}$ for which this property has not yet
  been established.  We consider two cases, according to whether $X$
  covers one or more elements of $\mc{Z}$.
  
  First assume $X$ covers only $Y$ in $\mc{Z}$.  Now $Y$ contains an
  independent set $J$ with $|J|=r(Y)$ by the choice of $X$.  Let $I$
  be the union of $J$ and any set of $r(X)-r(Y)$ elements of $X-Y$.
  If $I$ were dependent, it would contain a circuit $C$ with at most
  $r(X)$ elements. Thus, $r(\bar{C})<r(X)$, so, by
  (\ref{thm:axioms}.2), $\bar{C}\subsetneq X$.  Since $J$ is
  independent, $\bar{C}\not\subseteq Y$.  These conclusions contradict
  the assumption that $X$ covers only one set in $\mc{Z}$.  Thus, $I$
  is an independent subset of $X$ with $r(X)$ elements, as desired.
  
  Now assume $X$ covers $X_1$ and $X_2$ (and perhaps more sets) in
  $\mc{Z}$. Let $B'$ be a basis of $X_1\meet X_2$ and let $B''$ be
  $B'\cup \bigl((X_1\cap X_2)-(X_1\meet X_2)\bigr)$; by
  (\ref{thm:axioms}.3), $B''$ is a basis of $X_1\cap X_2$.
  Choose bases $B_1$ of $X_1$ and $B_2$ of $X_2$, both containing
  $B''$, and let $B$ consist of $B_1$ and any $r(X)-r(X_1)$ elements
  of $B_2-B''$.  (By condition (Z3), $|B_2-B''|\geq r(X)-r(X_1)$.)
  Thus, $|B|=r(X)$.  We claim that $B$ is independent.
  
  Assume, to the contrary, that $B$ contains a circuit $C$.  Since
  $C\not\subseteq X_1$, we have $X_1\subsetneq X_1\join
  \bar{C}\subseteq X$; since $X$ covers $X_1$, we have $X_1\join
  \bar{C}=X$.  Also, $r(X_1\meet \bar{C})=r_{\mc{Z}}(X_1\meet
  \bar{C})$ since $X_1\meet \bar{C} \subsetneq X$.  Applying
  property~(Z3) to $X_1$ and $\bar{C}$ and using (\ref{thm:axioms}.3)
  gives
  \begin{align}
  r(X_1)+|C|-1 \,&\geq \, r(X)+r(X_1\meet \bar{C}) +|(X_1\cap \bar{C})-
   (X_1\meet \bar{C})|\notag\\ 
   &= \, r(X)+r_{\mc{Z}}(X_1\meet \bar{C}) +|(X_1\cap \bar{C})-
   (X_1\meet \bar{C})|\notag\\ 
  & =\, r(X)+r_{\mc{Z}}(X_1\cap \bar{C})\notag\\
  &\geq \, r(X)+|X_1\cap C|\notag\\
  &\geq \, r(X)+ |C|-\bigl(r(X)-r(X_1)\bigr)\notag \\
   & = \, r(X_1)+|C|.\notag
  \end{align}
  (The last inequality follows by the construction of $B$.)  This
  contradiction implies that $B$ contains no circuit and hence is
  independent, as needed.
\end{proof}

The next two steps show that the sets in $\mc{Z}$ are cyclic flats of
$M_{\mc{Z}}$.

\begin{quote}
  {\bf (\ref{thm:axioms}.5)} Each set $X$ in $\mc{Z}$ is a flat of
  $M_{\mc{Z}}$.
\end{quote}

\begin{proof}[Proof of (\ref{thm:axioms}.5).]
  To show this, let $x$ be in a circuit $C$ of $M_\mc{Z}$ with
  $C\subseteq X\cup x$; we need to show that $x$ is in $X$.  Property
  (Z3) applied to $X$ and $\bar{C}$, together with
  (\ref{thm:axioms}.3) and (\ref{thm:axioms}.4), give the inequalities
  \begin{align}
    r(X) + |C|-1 \,& \geq \, r(X\join \bar{C}) + r(X\meet \bar{C}) +
    |(X\cap \bar{C}) - (X\meet \bar{C})|
    \notag\\  & = \, r(X\join \bar{C}) +r_{\mc{Z}}(X\cap \bar{C})\notag\\
    & \geq \, r(X\join \bar{C}) +|C|-1,\notag
  \end{align}
  where the last inequality follows since $C-x$ is an independent
  subset of $X\cap \bar{C}$. Hence $r(X)\geq r(X\join \bar{C})$, so
  $X=X\join \bar{C}$.  Thus $x\in X$ since $x\in \bar{C}$.
\end{proof}


\begin{quote}
  {\bf (\ref{thm:axioms}.6)} Each set $X$ in $\mc{Z}$ is a union of
  circuits and so is a cyclic flat.
\end{quote}

\begin{proof}[Proof of (\ref{thm:axioms}.6).]
  The claim holds for $0_\mc{Z}$ since, by property~(Z1), its elements
  are loops.  Let $X$ be a minimal set in $\mc{Z}$ for which the claim
  has not yet been verified and let $X$ cover the set $Y$ in $\mc{Z}$.
  Since $Y$ is a union of circuits, we need only show that each
  element $x$ in $X-Y$ is in a circuit that is contained in $X$.  Let
  $B$ be a basis of $Y$ and let $I$ be a set of $r(X)-r(Y)$ elements
  of $X-(Y\cup x)$, which exists by property~(Z2).  The subset $B\cup
  I\cup x$ of $X$ has $r(X)+1$ elements and so must contain a circuit
  $C$; we claim that $x$ is in $C$.  Assume, to the contrary, that $x$
  is not in $C$. Therefore $r(X)-r(Y)\geq |C-Y|$.  Since
  $\bar{C}\subseteq X$ and $C$ contains some elements of $I\cup x$, we
  have $Y\subsetneq Y\join \bar{C}\subseteq X$; that $X$ covers $Y$
  forces $Y\join \bar{C}=X$.  Applying property~(Z3) to $Y$ and
  $\bar{C}$, we obtain
  \begin{align}
    r(Y)+|C|-1\,& \geq \, r(X)+r(Y\meet \bar{C}) +|(Y\cap \bar{C})-
    (Y\meet \bar{C})|\notag\\
    & =\, r(X)+r_{\mc{Z}}(Y\cap \bar{C})\notag\\
    &\geq \,r(X)+|Y\cap C|.\notag
  \end{align}
  Combining this with the inequality $r(X)-r(Y)\geq |C-Y|$ noted above
  gives the contradiction $|C|-1\geq |C|$. Thus, $x$ is in $C$, so $X$
  is a union of circuits.
\end{proof}

Lastly, we show that the cyclic flats of $M_{\mc{Z}}$ are in $\mc{Z}$.


\begin{quote}
  {\bf (\ref{thm:axioms}.7)} Each cyclic flat $X$ of $M_\mc{Z}$ is in
  $\mc{Z}$.
\end{quote}

\begin{proof}[Proof of (\ref{thm:axioms}.7).]
  The assertion clearly holds for the least cyclic flat of $M_\mc{Z}$.
  Let $X$ be a minimal cyclic flat of $M_\mc{Z}$ for which the claim
  has not yet been verified and assume that $X$ covers $Y$ in
  $\mc{Z}(M_\mc{Z})$.  Thus, $Y$ is in $\mc{Z}$.  Let $B$ be a basis
  of $Y$ and let $T$ be any set of $r_{\mc{Z}}(X)-r_{\mc{Z}}(Y)+1$
  elements of $X-Y$.  Since $X$ is cyclic and covers $Y$, it follows
  that $B\cup T$ contains a circuit $C$ and that $C\cap (X-Y)$ is $T$.
  Note that $X=\cl_{M_{\mc{Z}}}(Y\cup C)$ since $X$ covers $Y$
  and contains $C$.  Since $Y$ and $\bar{C}$ are in $\mc{Z}$, their
  join $Y\join \bar{C}$ in $\mc{Z}$ is a cyclic flat of $M_{\mc{Z}}$
  by (\ref{thm:axioms}.6), and this flat must contain the smallest
  flat that contains $Y$ and $C$, which is $X$.  We will show that
  $Y\join \bar{C}$ is $X$, which implies that $X$ is in $\mc{Z}$, as
  claimed.  It suffices to prove the inequality $r(\bar{C}\join Y)
  \leq r_\mc{Z}(X)$.  Note that in the semimodular inequality in
  $M_\mc{Z}$,
  $$r_\mc{Z}(\bar{C}\cup Y) + r_\mc{Z}(\bar{C}\cap Y) \leq
  r_\mc{Z}(\bar{C}) + r_\mc{Z}(Y),$$
  the right side is $|C|-1+r(Y)$,
  that is, $\bigl(|C\cap Y| + r_\mc{Z}(X) - r(Y)\bigr) + r(Y)$, which
  is $|C\cap Y| + r_\mc{Z}(X)$.  The left side is clearly at least as
  large, so we have
  \begin{equation}\label{s7}
  r_\mc{Z}(\bar{C}\cup Y) + r_\mc{Z}(\bar{C}\cap Y) =
  r_\mc{Z}(\bar{C}) + r_\mc{Z}(Y).
  \end{equation}
  By (\ref{thm:axioms}.3) and property~(Z3), we have
  $$r(\bar{C}\join Y) +r_{\mc{Z}}(\bar{C}\cap Y)= r(\bar{C}\join Y) +
  r(\bar{C}\meet Y) + |(\bar{C}\cap Y) - (\bar{C}\meet Y)| \leq
  r(\bar{C})+r(Y),$$
  which, with Eqn.~(\ref{s7}) and statement
  (\ref{thm:axioms}.4), gives $r(\bar{C}\join Y) \leq
  r_\mc{Z}(\bar{C}\cup Y)$.  Now $ r_\mc{Z}(C)= r_\mc{Z}(\bar{C})$, so
  $r_{\mc{Z}}(X)= r_\mc{Z}(C\cup Y)= r_\mc{Z}(\bar{C}\cup Y)$, so
  the desired inequality $r(\bar{C}\join Y) \leq r_\mc{Z}(X)$ follows.
\end{proof}

Thus, as needed, from $\mc{Z}$ and $r:\mc{Z}\rightarrow \mathbb{Z}$,
we have constructed a matroid $M_\mc{Z}$ for which $\mc{Z}(M_\mc{Z})$
is $\mc{Z}$ and $r$ is the restriction of the rank function of
$M_\mc{Z}$ to $\mc{Z}$.
\end{proof}

In Theorem~\ref{thm:axioms}, we rediscovered a result of
Sims~\cite[Chapter~3, Theorem~2.2]{juliethesis}.  Indeed, she proves
the theorem for rank-finite independence spaces in the case that $S$
is in $\mc{Z}$. (Our proof also applies in that context; both
approaches also apply when there are only finitely many elements in
$S$ that are not in the greatest set in $\mc{Z}$.)  Her proof consists
of defining a function $r'$ on the subsets of $S$ by
$$r'(A)=\min\{r(F)+|A-F|\,:\,F\in\mc{Z}\}$$
and then proving (i) $r'$
satisfies the rank axioms of a rank-finite independence space
(specifically, $r'(\emptyset)=0$, the unit increase property,
semimodularity, and that each subset $A$ of $S$ has a finite subset
$B$ with $r'(A)=r'(B)$), (ii) $r(F)=r'(F)$ for $F\in \mc{Z}$, and
(iii) the sets in $\mc{Z}$ are precisely the cyclic flats of the
resulting independence space.  One intermediate result that proves
quite useful is that $r(G)<r(F)+|G-F|$ for all pairs of distinct sets
$F$ and $G$ in $\mc{Z}$.

We close this section by noting that viewing matroid operations from
the perspective of cyclic flats can suggest extensions of familiar
operations.  For instance, relaxing a circuit-hyperplane $H$ consists
of removing $H$ from the lattice of cyclic flats.  It follows from
Theorem~\ref{thm:axioms} that one can remove any cyclic flat that,
like a circuit-hyperplane, is comparable only to the least and
greatest cyclic flats.


\section{The Free Product from the Perspective of Cyclic 
Flats}\label{sec:free}


Crapo and Schmitt~\cite{free} introduced a noncommutative matroid
operation called the free product.  Many interesting properties and
applications are treated in~\cite{free,free2}.  In this section we
show that cyclic flats give a transparent path to the free product.
We connect this approach with the definition in~\cite{free}, which
uses independent sets, and we use a corollary of this link to compute
Tutte polynomials of free products.

Our starting point for the free product is the description of the
cyclic flats given in \cite[Proposition~6.1]{free2}.  The new
observation is that from this description of the cyclic flats, it is
obvious that properties (Z0)--(Z3) of Theorem~\ref{thm:axioms} hold.
Thus, cyclic flats give a simple way to define the free product.

\begin{thm}\label{thm:fp}
  Let $M$ and $N$ be matroids on disjoint ground sets.  Set
  $$\mc{Z}'= \{X \,:\, X\in\mc{Z}(M),\, X\ne E(M)\} \cup \{E(M)\cup Y
  \,:\, Y\in\mc{Z}(N),\, Y\ne\emptyset\}.$$
  Let $\mc{Z}$ be
  $\mc{Z}'\cup \{E(M)\}$ if $M$ has no isthmuses and $N$ has no loops;
  otherwise let $\mc{Z}$ be $\mc{Z}'$.  Set $r(X)=r_M(X)$ and
  $r(E(M)\cup Y)=r(M)+r_N(Y)$.  The pair $(\mc{Z},r)$ satisfies
  conditions \emph{(Z0)--(Z3)} of Theorem~\ref{thm:axioms} and so
  defines a matroid on $E(M)\cup E(N)$.
\end{thm}

We take this as the definition of the free product $M\frp N$.  Note
that if $M$ has an isthmus and $N$ has a loop, then $\mc{Z}$ is
isomorphic to the linear sum of $\mc{Z}(M)$ and $\mc{Z}(N)$; otherwise
the difference between $\mc{Z}$ and the linear sum is that the
greatest element of $\mc{Z}(M)$ and the least element of $\mc{Z}(N)$
are, in effect, identified in $\mc{Z}$.  Only property (Z3) requires
even a minor observation: the only nontrivial instances of this
property are either instances of this property in $M$ or shifts by
$2\,r(M)$ of instances of this property in $N$.

Two cases of the free product, treated in~\cite{free,free2}, deserve
special mention.  The lattice of cyclic flats of $M\frp U_{0,1}$
agrees with that of $M$ except that $E(M)\cup E(U_{0,1})$ is a cyclic
flat in $M\frp U_{0,1}$ while $E(M)$ is not; ranks are unchanged; this
is the free extension of $M$.  The lattice of cyclic flats of
$U_{1,1}\frp M$ is formed from that of $M$ by augmenting every
nonempty set by $E(U_{1,1})$; ranks increase by $1$; this is the free
coextension of $M$.  Theorem~\ref{thm:fp} and this view of free
coextension give the following informal geometric description of the
free product: $M\frp N$ is formed by taking $r(M)$ free coextensions
of $N$, say by $e_1,e_2,\ldots,e_{r(M)}$, gluing $M$ freely into the
flat spanned by $e_1,e_2,\ldots,e_{r(M)}$, and then deleting
$e_1,e_2,\ldots,e_{r(M)}$.

We now show that the independent sets of $M\frp N$, which were used to
define this operation in~\cite[Proposition~1]{free}, can be identified
easily from the cyclic flat perspective.  We let $\nu(X)$ denote the
nullity, $|X|-r(X)$, of $X$.

\begin{thm}\label{thm:ind} Let $\mc{I}_M$ be the collection of
  independent sets of $M$.  The collection of independent sets of
  $M\frp N$ is given by
  $$\{X\cup Y \,:\,X \in \mc{I}_M,\, Y\subseteq E(N), \text{ and \ }
  \nu_N(Y)\leq r(M)-|X| \}.$$
\end{thm}

\begin{proof}
  By Lemma~\ref{lem:cyc+rank}, for $X\subseteq E(M)$ and $Y\subseteq
  E(N)$ the set $X\cup Y$ is independent in $M\frp N$ if and only if
  $|(X\cup Y)\cap Z|\leq r(Z)$ for every cyclic flat $Z$ of $M\frp N$.
  This condition, with $Z$ ranging over $\mc{Z}(M)$ and (if it is in
  $\mc{Z}(M\frp N)$) $E(M)\cup \cl_N(\emptyset)$, is equivalent to $X$
  being independent in $M$.  The other cyclic flats of $M\frp N$ have
  the form $E(M)\cup Z'$ where $Z'$ is in $\mc{Z}(N)$. Since
  $r(E(M)\cup Z')=r(M)+r_N(Z')$, the inequality of interest is
  equivalent to $|X|+ |Y\cap Z'|\leq r(M)+r_N(Z')$, or $|Y\cap
  Z'|-r_N(Z') \leq r(M) - |X|$.  The proof is completed by noting that
  the maximum of $|Y\cap Z'|-r_N(Z')$ over all cyclic flats of $N$ is
  $\nu_N(Y)$.  Indeed, we have $|Y\cap Z'|-r_N(Z')\leq \nu_N(Y\cap
  Z')\leq \nu_N(Y)$, and if $Z'$ is the largest cyclic flat of $N$
  contained in $\cl_N(Y)$, then $|Y\cap Z'|-r_N(Z')$ is $\nu_N(Y)$.
\end{proof}

By~\cite[Proposition~7.2]{free2}, the free product $M\frp N$ is the
freest matroid (in the weak order) on $E(M)\cup E(N)$ whose
restriction to $E(M)$ is $M$ and whose contraction to $E(N)$ is $N$;
also, the direct sum $M\oplus N$ is the least matroid in the weak
order with these specified minors.  Both operations are simple from
the perspective of cyclic flats: the lattice of cyclic flats of
$M\oplus N$ is isomorphic to the direct product of the lattices of
cyclic flats of $M$ and $N$ and, as noted above, the lattice of cyclic
flats of $M\frp N$ is related to the linear sum of those for $M$ and
$N$.

Recall that the Tutte polynomial of a matroid $M$ is given by
$$t(M;x,y) = \sum_{A \subseteq E(M)}
(x-1)^{r(M)-r(A)}(y-1)^{\nu(A)}.$$
We now show that, given the Tutte
polynomials of two matroids, it is easy to compute the Tutte
polynomial of their free product.  However, unlike the simple formula
$t(M\oplus N;x,y)=t(M;x,y)t(N;x,y)$ for direct sums, for the free
product we get a formula for each coefficient in the Tutte polynomial.
We use the next lemma~\cite[Proposition 3.5]{free}, which follows from
Theorem~\ref{thm:ind}.

\begin{lemma}\label{lem:fprank}
  For $X\subseteq E(M)$ and $Y\subseteq E(N)$, we have
  $$r_{M\frp N}(X\cup
  Y)=r_M(X)+r_N(Y)+\min\{r(M)-r_M(X),\nu_N(Y)\}.$$
\end{lemma} 

\begin{thm}\label{poly}
  The Tutte polynomial of a free product can be computed in polynomial
  time in the size of its ground set from the Tutte polynomials of its
  factors.
\end{thm}

\begin{proof}
  First note that $r(M\frp N)=r(M)+r(N)$, so for $X\subseteq E(M)$ and
  $Y\subseteq E(N)$, Lemma~\ref{lem:fprank} gives
  $$r(M\frp N)-r_{M\frp N}(X\cup Y)=r(M)-r_M(X)+r(N)-r_N(Y)
  -\min\{r(M)-r_M(X),\nu_N(Y)\}$$
  and
  $$\nu_{M\frp N}(X\cup Y)=\nu_M(X)+\nu_N(Y)-
  \min\{r(M)-r_M(X),\nu_N(Y)\}.$$
  We prove the theorem using the Whitney rank generating function
  $R(M;x,y)$, which is $t(M;x+1,y+1)$.  Thus, the coefficient of
  $x^iy^j$ in $R(M;x,y)$ is the number of subsets $A$ of $E(M)$ with
  $r(M)-r_M(A)=i$ and $\nu_M(A)=j$.  Let $R(M;x,y)$ and $R(N;x,y)$ be
  $$
  R(M;x,y)=\sum\limits_{\substack{ 0\leq i \leq r(M) \\
      0\leq j \leq \nu(M)}} a_{ij}x^iy^j
  \qquad
  \textrm{ and }
  \qquad
  R(N;x,y)=  \sum \limits_{\substack{ 0\leq k \leq r(N) \\
      0\leq l \leq \nu(N)}} b_{kl}x^ky^l.
  $$
  The number of coefficients in these polynomials is at most
  $(\nu(M)+1)(r(M)+1)$ and $(\nu(N))+1)(r(N)+1)$, respectively.  By
  the equations in the first sentence, the coefficient of $x^p y^q$ in
  $R(M\frp N;x,y)$ is
  $$\sum\limits_{\substack{ i,j,k,l: \\
      i+k-\min\{i,l\} = p \\ j+l-\min\{i,l\} = q }} a_{ij}b_{kl}.$$
  Thus, each of the $(|E(M)|+|E(N)|+1)^2$ or fewer coefficients of
  $R(M\frp N;x,y)$ can be computed with polynomially many steps from
  $R(M;x,y)$ and $R(N;x,y)$, so $R(M\frp N;x,y)$ can be computed from
  $R(M;x,y)$ and $R(N;x,y)$ in polynomial time in $|E(M)|+|E(N)|$.
\end{proof}


\section{Cyclic Width}\label{sec:cw}


Order-theoretic properties of the lattice of cyclic flats can be used
to define some classes of matroids; this section begins to explore
some such classes.

The width of a lattice is the maximal cardinality of an antichain
(i.e., a set of incomparable elements) in the lattice, so we define
the \emph{cyclic width} of a matroid $M$ to be the width of
$\mc{Z}(M)$.  For instance, the matroids in Figure~\ref{fig:lattice}
have cyclic width~4.  We use $\cw{k}$ to denote the class of all
matroids whose cyclic width is at most $k$.  By
Lemma~\ref{lem:nested}, the class of nested matroids is $\cw{1}$.  We
first show that $\cw{k}$ is closed under several basic matroid
operations.

\begin{thm}\label{thm:cwk}
  The class $\cw{k}$ is closed under duals, minors, and free products.
\end{thm}

\begin{proof}
  Closure under duals and free products follows from the views of
  these operations in terms of cyclic flats given in
  Sections~\ref{sec:intro} and~\ref{sec:free}.  For the result on
  minors, it suffices to show that if $M\in\cw{k}$ and $x\in E(M)$,
  then $M\del x$ is in $\cw{k}$.  If this were false, then there would
  be flats $X_1,X_2,\ldots,X_{k+1}$ in $\mc{Z}(M\del x)$ with
  $X_i\not\subseteq X_j$ for $i\ne j$.  Note that exactly one of $X_i$
  and $X_i\cup x$ would be in $\mc{Z}(M)$; let
  $X'_1,X'_2,\ldots,X'_{k+1}$ be these flats in $\mc{Z}(M)$.  That $M$
  is in $\cw{k}$ gives the inclusion $X'_i\subseteq X'_j$ for some
  distinct $i$ and $j$.  This inclusion gives the contradiction
  $X_i\subseteq X_j$, so $M\del x$ is, as claimed, in $\cw{k}$.
\end{proof}

Recall that the truncation $T(M)$ of $M$ is $(M+e)/e$; this is the
matroid on $E(M)$ whose bases are the independent sets of $M$ of size
$r(M)-1$.  The Higgs lift is the dual operation.

\begin{cor}
  The class $\cw{k}$ is closed under free extension, free coextension,
  truncation, and the Higgs lift.
\end{cor}

Note that $\cw{k}$ is not closed under direct sums.

We make some remarks about the excluded minors for $\cw{k}$.  The
excluded minors for $\cw{1}$ were shown in~\cite{opr} to be $P_n$, for
$n\geq 2$, where $P_n$ is the (iterated) truncation to rank $n$ of
$U_{n-1,n}\oplus U_{n-1,n}$.  Likewise, $\cw{k}$ has infinitely many
excluded minors, including, for $n\geq 2$, the truncation to rank $n$
of the direct sum of $k+1$ copies of $U_{n-1,n}$; the dual is also an
excluded minor.  For $k>1$, there are other excluded minors for
$\cw{k}$, such as the truncation to rank~$3$ of the $(k+1)$-whirl.
Not all excluded minors for $\cw{k}$ have cyclic width $k+1$.  Indeed,
as $k$ grows, the difference between $k$ and the cyclic width of an
excluded minor for $\cw{k}$ can be arbitrarily large, as the following
example shows.  Let $M_n$ be the truncation to rank three of the
rank-$(n+1)$ binary projective geometry $\PG(n,2)$.  We claim that
$M_n$ is an excluded minor for $\cw{k}$, where $k$ is
$(2^{n+1}-4)(2^n-1)/3$, and the cyclic width of $M_n$ exceeds $k$ by
$2^n-1$.  The cyclic width of $M_n$ is the number of lines of
$\PG(n,2)$, which is $(2^{n+1}-1)(2^n-1)/3$.  Every element of
$\PG(n,2)$ is on $2^n-1$ lines, all of which have three points, so
every single-element deletion of $M_n$ has cyclic width
$(2^{n+1}-1)(2^n-1)/3 - (2^n-1)$ (which is $k$) and every
single-element contraction of $M_n$ has cyclic width $2^n-1$.

We now turn to $\cw{1}$.  The class of nested matroids has many
interesting properties, which partly explains why these matroids have
been introduced a number of times in different contexts
(see~\cite[Section~4]{lpm2}).  Recall that a class of matroids is
well-quasi-ordered if it contains no infinite antichain in the minor
order, that is, there is no infinite set of matroids in the class none
of which is isomorphic to a minor of another.  Theorem~\ref{thm:wqo}
shows that $\cw{1}$ is an example (apparently the first known) of a
well-quasi-ordered class of matroids that, as mentioned above, has
infinitely many excluded minors.  This result, along with well-known
examples, shows that there is no connection between the following two
properties that a minor-closed class $\mc{M}$ of matroids may have:
(a) $\mc{M}$ is well-quasi-ordered; (b) $\mc{M}$ has a finite set of
excluded minors.

To prove that nested matroids are well-quasi-ordered, we use Higman's
theorem \cite{higman}, which is stated in Lemma~\ref{lem:Nt} (see
also, e.g., \cite[Theorem~5.2]{thom}).  Recall that if $X$ is
quasi-ordered, then the set of finite sequences in $X$ is
quasi-ordered as follows: for $\mathbf{x}=x_1,x_2,\ldots,x_m$ and
$\mathbf{y}=y_1,y_2,\ldots,y_n$, set $\mathbf{x}\leq \mathbf{y}$ if
there are integers $1\leq i_1 <i_2 <\cdots< i_m\leq n$ with $x_j\leq
y_{i_j}$ for $1\leq j\leq m$.

\begin{lemma}\label{lem:Nt}
  If $X$ is well-quasi-ordered, then so is the set of finite sequences
  of $X$.
\end{lemma}

\begin{thm}\label{thm:wqo}
  Nested matroids are well-quasi-ordered.
\end{thm}

\begin{proof}
  We treat isomorphic matroids as equal.  A nested matroid $M$ is
  formed from the empty matroid by applying two operations: (i) adding
  an isthmus and (ii) adding an element freely.  Let $\mathbf{v}_M$ be
  the corresponding sequence of $i$'s and $f$'s.  Consider the order
  on $\{i,f\}$ in which $i$ and $f$ are incomparable; this is a
  well-quasi-order.  The theorem follows from Lemma~\ref{lem:Nt} by
  showing that for nested matroids $M$ and $N$, if $\mathbf{v}_N\leq
  \mathbf{v}_M$, then $N$ is a minor of $M$.  Indeed,
  $\mathbf{v}_N\leq \mathbf{v}_M$ means that $\mathbf{v}_N$ is a
  subsequence of $\mathbf{v}_M$, so to obtain $N$, from $M$ remove the
  elements that do not contribute to this subsequence, deleting those
  that were added freely and contracting those that were added as
  isthmuses.
\end{proof}

A.~M.~H.~Gerards~\cite{bert} has noted that Theorem~\ref{thm:wqo} is a
concrete instance of a general result: any quasi-order $X$ that is not
a well-quasi-order and that has no infinite descending chains contains
a well-quasi-order $Y$ for which infinitely many elements in $X$ are
minimal in $X-Y$.  This follows by letting $Y$ be the set of all
elements of $X$ that are smaller than at least one element of a fixed
minimal bad sequence (as defined in~\cite{diestel}). 

Note that $\cw{2}$ is not well-quasi-ordered; the excluded minors for
$\cw{1}$ (the matroids $P_n$ mentioned above) are an infinite
antichain in $\cw{2}$.

It appears that few matroids in $\cw{2}$ have been studied previously.
Acketa~\cite{acketa} proved that any matroid whose lattice of cyclic
flats is a product of two $2$-element chains is transversal.  Our next
result shows that this conclusion holds for all matroids in $\cw{2}$.
This result is a consequence of the following characterization of
transversal matroids due to Ingleton~\cite{ingleton}, which refines a
result of Mason.

\begin{lemma}\label{lem:mason-ingleton}
  A matroid is transversal if and only if for every nonempty family
  $(X_1,\ldots,X_n)$ of cyclic flats of $M$,
\begin{equation}\label{mi}
r(X_1\cap X_2\cap \cdots \cap X_n)\leq \sum_{J\subseteq \{1,2,\ldots,n\}}
  (-1)^{|J|+1} r\Bigl(\bigcup_{j\in J} X_j \Bigr).
\end{equation}
\end{lemma} 

In Lemma~\ref{lem:mason-ingleton}, it suffices to consider antichains
of cyclic flats since if $X_i\subseteq X_j$, then omitting $X_j$ does
not change either side of inequality~(\ref{mi}).  Indeed, the terms on
the right side that include $X_j$ cancel via the involution that
adjoins or omits $X_i$.

\begin{thm}\label{thm:cw2}
  Matroids in $\cw{2}$ are bitransversal.
\end{thm}

\begin{proof}
  By duality it is enough to show that matroids in $\cw{2}$ satisfy
  the condition in Lemma~\ref{lem:mason-ingleton}.  In $\cw{2}$,
  antichains of cyclic flats contain at most two flats.  Equality
  holds in inequality~(\ref{mi}) for $n=1$. For $n=2$,
  inequality~(\ref{mi}) is the semimodular inequality, $r(X_1\cap
  X_2)\leq r(X_1)+r(X_2)-r(X_1\cup X_2)$.
\end{proof}

It follows from the definition of nested matroids that there are, up
to isomorphism, $2^n$ nested matroid on $n$ elements.  In contrast,
the class $\cw{2}$ is superexponential, as shown by the following
result due to O.~Gim\'enez~\cite{superexponential}.

\begin{thm}\label{thm:omer}
In $\cw{2}$, there are at least $n!$ matroids on $4n+5$ elements that all
have isomorphic lattices of cyclic flats. 
\end{thm}

\begin{proof}
  Let $S$ be $\{a,a',a'',b,b',x_1,\ldots,
  x_n,y_1,\ldots,y_n,z_1,\ldots, z_n, w_1,\ldots, w_n \}$.  For a
  permutation $\sigma$ of $\{1,2,\ldots,n\}$, let $\mc{Z}_{\sigma}$ be
  $\{\emptyset, A_0, A_1,\ldots,A_n, B_0, B_1,\ldots, B_n, S\}$ where
  $$
  \begin{array}{l} A_0=\{a,a',a'',x_1, \ldots, x_n\}, \\
    B_0=\{b,b',y_1,\ldots, y_n\},\\
    A_i=A_{i-1}\cup \{z_i,w_i\}, \quad \mathrm{ for\ } 1\leq i \leq n,\\
    B_i=B_{i-1}\cup \{z_i,w_{\sigma(i)}\}, \quad \mathrm{ for\ } 1\leq
    i \leq n.\\
  \end{array}$$
  Note that the only inclusions among these sets are $\emptyset\subset
  A_0\subset \cdots \subset A_n \subset S$ and $\emptyset\subset
  B_0\subset \cdots \subset B_n \subset S$. Also, $\mc{Z}_\sigma$ is a
  lattice.  Define $r:\mc{Z}_{\sigma}\rightarrow \mathbb{Z}$ by
  $r(\emptyset)=0$, $r(A_i)=r(B_i)=n+1+i$ for $0\leq i \leq n$, and
  $r(S)=2n+2$. We show that the pair $(\mc{Z}_{\sigma},r)$ satisfies
  conditions (Z1)--(Z3) of Theorem~\ref{thm:axioms} and so defines a
  matroid $M_{\sigma}$.  Conditions (Z1) and (Z2) hold by
  construction.  For condition (Z3), we need
  $$r(A_i) +r(B_j)\geq r(A_i\join B_j) +r(A_i\meet B_j) +|(A_i\cap
  B_j)- (A_i \meet B_j)|$$
  for all $i,j$ between $0$ and $n$.  The
  left side is $n+i+1+n+j+1$; since $A_i\join B_j=S$ and $A_i\meet B_j=
  \emptyset$, the right side is at most $2n+2+2\min\{i,j\}$, so the
  required inequality holds.  Finally, note that $\sigma$ can be
  recovered from $\mc{Z}(M)$ for any matroid $M$ isomorphic to
  $M_\sigma$, so different permutations give nonisomorphic matroids.
\end{proof}

The construction in this proof can be adapted (using $k-1$
permutations) to show that there are at least $(n!)^{k-1}$
nonisomorphic matroids in $\cw{k}$ on a set of
$(k+2)n+\binom{k+2}{2}-1$ elements, all with isomorphic lattices of
cyclic flats.

We now consider the matroids in $\cw{3}$ that are binary, ternary, or
have rank~$3$.  (See~\cite{ingleton} or~\cite{oxley} for the
definitions of gammoids and base-orderable matroids.)

\begin{thm}\label{thm:cw3}
  Binary matroids in $\cw{3}$ are direct sums of series-parallel
  networks.  Rank-$3$ matroids in $\cw{3}$ are base-orderable.
  Ternary matroids in $\cw{3}$ are gammoids.
\end{thm}

\begin{proof}
  Since $M(K_4)$ has cyclic width~$4$ and $\cw{3}$ is minor-closed, no
  matroid in $\cw{3}$ has an $M(K_4)$-minor.  The first and second
  assertions follow since (i) binary matroids with no $M(K_4)$-minor
  are direct sums of series-parallel networks~\cite[Theorem
  13.4.9]{oxley} and (ii) a rank-$3$ matroid is base-orderable if and
  only if it has no $M(K_4)$-restriction~\cite[Theorem 14]{ingleton}.
  The third statement follows similarly since the ternary excluded
  minors of ternary gammoids have cyclic width at least~4 (the
  excluded minors for ternary gammoids were found in~\cite{tern_gam}).
\end{proof}


\section{Open Problems}


We close with several problems suggested by the topics of this paper.

Theorem~\ref{thm:alllattices} states that every lattice is isomorphic
to the lattice of cyclic flats of a matroid.  The problem is more
subtle if we drop the phrase ``isomorphic to''.

\begin{prob}
  Characterize the lattices of sets that are lattices of cyclic flats
  of matroids.
\end{prob}

In other words: for which lattices $\mc{Z}$ of sets can one find a
function $r:\mc{Z}\rightarrow\mathbb{Z}$ so that conditions (Z1)--(Z3)
of Theorem~\ref{thm:axioms} hold?

\begin{prob}
  What is the largest number of cyclic flats that a matroid on $n$
  elements can have?
\end{prob}

The direct sum of $\lfloor n/2\rfloor$ copies of $U_{1,2}$, along with
$U_{1,1}$ if $n$ is odd, has $n$ elements and $2^{\lfloor n/2
  \rfloor}$ cyclic flats.  For $n\geq 7$, the number of cyclic flats
is greater in binary spikes (with or without a tip, according to the
parity of $n$), and, for sufficiently large $n$, even greater in
certain ternary matroids that are similar to spikes.  We are not
currently aware of any plausible candidates for matroids that may
maximize the numbers of cyclic flats.

Another line of research is to investigate which properties of
$\cw{1}$ have counterparts for other classes $\cw{k}$.  For instance,
matroids in $\cw{1}$ and $\cw{2}$ are transversal, and those in
$\cw{1}$ have presentations by particularly simple set systems
(see~\cite{lpm1,opr}).  Is there an interesting description of some
presentations of matroids in $\cw{2}$?  Under certain linear orders,
the broken circuit complex of a matroid in $\cw{1}$ is the
independence complex of another matroid in $\cw{1}$ (see~\cite{lpm1}).
Do matroids in $\cw{2}$ have a similar property?  Are all matroids in
$\cw{3}$ gammoids?  A weaker question is this: are all matroids in
$\cw{3}$ base-orderable?  (This is not the case for $\cw{4}$ since
$M(K_4)$ is both the smallest matroid that is not a gammoid and the
smallest matroid that is not base-orderable.)  

\bigskip

\begin{center}
{\sc Acknowledgements}
\end{center}

Bert Gerards and Omer Gim\'enez allowed us to include their
observations.  Joseph Kung provided a number of useful comments,
especially related to Dilworth's embedding theorem.  Julie Sims kindly
scanned her thesis so that we could sketch her proof of
Theorem~\ref{thm:axioms}.  We are very grateful to them for their
assistance.

\end{document}